# Optimization of Quadratic Sieve Algorithm Implementation for Large Integer Factorization

ZIHAN GUAN[1], XIAODONG ZHUANG[1], NIKOS E. MASTORAKIS[2]
[1]Electronics Information College,
Qingdao University, Qingdao,
CHINA

[2]Technical University of Sofia,
Sofia,
BULGARIA

*Abstract:* - The quadratic sieve method is a core tool in number theory. In this paper, we present two optimization methods for the Quadratic Sieve algorithm. In the sieving process, the original polynomial root value accumulation step is changed from the original $d_i$ to the $md_i$ (m is a small integer), which can change the complexity from $O(n)$ to $O\left(\frac{n}{m}\right)$. Another optimization method is that for all parameter lookup steps, the original traversal lookup can be changed to an efficient binary search, which can change the complexity of the loop from $O(n)$ to $O(logn)$. This enhancement reduces the computational complexity of RSA modulus factorization in practical settings.

*Key-Words:* - integer factorization, quadratic sieve, parameter lookup tables, factor base, algorithm optimization, prime numbers, large integers.

## 1 Introduction

RSA is one of the most widely used cryptographic systems [1]. It is now used in many electronic systems, such as securing web traffic, e-mail, and some wireless devices [2]. Its security mainly depends on the difficulty of large integer factorization. In 1978, the RSA public-key cryptography system was proposed [3], the following decades, a variety of factorization [4] factor algorithms appeared, such as Pollard's rho method [5], elliptic curve [6] method, number field sieve [7] and quadratic sieve [8] method.

In this paper, the current quadratic sieve algorithm is analyzed, and its implementation is optimized. In the implementation, the original traversal algorithm is improved by step-size optimization and binary search. The algorithm optimization method adopted in this paper has greatly improved the efficiency of factoring large integers.

## 2 Brief Introduction of Quadratic Sieve

The basic quadratic sieve algorithm has origins which date back to 1980s [9], whose innovative factorization concepts provided its direct precursor. The algorithm decomposes large integers primarily through the following steps: constructing a quadratic congruence, establishing a factor base, performing the sieving process, solving a system of linear equations, and ultimately deriving factors. For large integers ranging from approximately 50 to 100 decimal digits, the Quadratic Sieve remains a practically valuable factorization method.

The quadratic sieve, a factoring algorithm [10] for composite integers $\mathrm{n}$, exploits properties of $\mathrm{Z}/\mathrm{nZ}$ exponentiation. It aims to find nontrivial integer solutions satisfying $X^2 \equiv Y^2 (\mathrm{mod}\ \mathrm{N})$ [11], thereby obtaining nontrivial factors of N by computing $gcd(\mathrm{x}-\mathrm{y},\mathrm{N})$ and $gcd(\mathrm{x}+\mathrm{y},\mathrm{N})$.

## 3 Analysis of Algorithm Implementation

Based on the above theory, this section mainly introduces the design and implementation of the Quadratic Sieve algorithm in detail. First, the input large integer undergoes preprocessing, and its optimal k value is determined using empirical formulas. The factor base is then constructed. The parameter configuration and polynomial construction are carried out. Following this, the smooth relations are obtained and the matrix is constructed for linear combination. Finally, factorization is achieved by computing the greatest common divisor (GCD). The programming implementation has four major parts. The implementation process is shown in Figure 1 (Appendix).

### 3.1 Data Pre-Processing (Small Factor Factorization and Perfect Power Judgment)

For the input integer $N$, parse the input digits into binary format. If $N$ has fewer than 64 bits, the factorization will be performed directly, and the factorization result will be returned. Otherwise, test whether N is a perfect power. If $N$ has the form of $x^y$, where x and y are integers, it can be directly decomposed. If $N$ does not have this form, it is further tested to see whether $N$ is prime. If it is composite, the quadratic sieve method is used to further decompose it.

#### (1) Module value $k$ and parameter configuration

To improve the smoothness probability and handle the case where the number of input bits is insufficient, the input integer $N$ is multiplied by the small prime $k$. According to the formula (1), the highest score k is selected.

$$\mathrm{S}(k) = \sum_{p\in \mathrm{B}} \frac{\log p}{p-1} \cdot \delta\left(\left(\frac{kN}{p}\right) = 1\right) \tag{1}$$

A parameter table is provided for different magnitudes, including factor base size, sieve interval length [12], and large prime threshold. According to the comparison between the bits number of $kN$ and the existing parameters, the interval of the parameter list is obtained by traversing. Linear interpolation within this interval establishes the relationship between $kN's$ bit-length and each parameter. Finally, the corresponding parameter value can be obtained according to the bits number of $kN$, which leads to the optimal parameter configuration.

#### (2) Construct the polynomials

The polynomial is constructed for the efficient multiple sieve.

$$f(x) = Ax^2 + 2Bx + C \tag{2}$$

To calculate $A$, first calculate the hypercube parameter $D$.

$$D = \frac{\sqrt{2kN}}{M} \tag{3}$$

$A$ is formed as a product of a subset of primes from the factor base, chosen to be closest to $D$. Let $B$ satisfy: $B^2 \equiv kN \bmod A$. By the Chinese remainder theorem, the $P_i$ of all the prime numbers that make up $A$ are substituted into the formula (4).

$$B_i \equiv \bmod P_i \tag{4}$$

Substitute the above data into the formula (5) to obtain $B$.

$$B = \sum_{i=1}^{s} B_i \cdot \left(\frac{A}{P_i}\right)\left(\left(\frac{A}{P_i}\right)^{-1} \bmod P_i\right) \bmod A \tag{5}$$

C is then derived as:

$$C = \frac{B^2 - kN}{A} \tag{6}$$

According to the above, the polynomial $f(x)$ is constructed. Then calculate the polynomial (7), and solve: $root_{1,2} \equiv \frac{-B \pm \sqrt{kN}}{A} \bmod p$.

$$f(x) \equiv 0 \bmod p \tag{7}$$

Keep $A$ unchanged but adjust $B$ to construct multiple polynomials, [13]. There is large amount of computation to find the root of the newly constructed polynomial. The root-value-update method is used to obtain these roots.

#### (3) The sieving operation

The smooth numbers of polynomial $f(x)$ can be selected by the sieve method, which greatly reduces the computation. First, the initialization settings are pre-processed, including the determination of the sieve interval $x \in \left[-\frac{M}{2}, \frac{M}{2}\right]$ and a sieving array and a flag array. Then we start from the root of the polynomial $f(x)$ and compute remainder when divided by primes in the factor base. The remainder is used as the step size $d_i$ in the following sieving process. And add the step size to the root positions. This can be implemented by marking the corresponding position of $x$ in the flag array. Then find the $x$ that is most likely to be factored by the primes in the factor base. Calculate the corresponding $f(x)$ and decompose it in the factor base. The exponent of each decomposed prime is recorded in the matrix named as $power1$. Next, replace the odd number in the matrix $power1$ with 1 and the even number with 0 to construct matrix $power2$. Find the rows of the matrix, and add them in order to get a zero vector, which are used to construct the square congruence formula (8).

$$x^2 \equiv y^2 \bmod N \tag{8}$$

Finally, at least one of $\gcd(x \pm y, N)$ is a non-trivial factor, and the factorization is completed.

# 4 Optimization of Programming Implementation

In this paper, the optimization is mainly for the sieving process. In the original program, root values $root_{1,2}$ of the polynomial $f(x)$ accumulate by step size $d_i$. The accumulation is implemented step by step, which is inefficient and has weak memory access continuity. It affects the overall efficiency of factorization. After the optimization, the step size for processing the remaining value position is changed to $md_i$ ($m$ is a small integer, here $m$ is 4 in the optimized program). The number of sieve cycles is reduced by 75%, and the time complexity of the cycle decreases from $O(n)$ to $O\left(\frac{n}{4}\right)$.

To better illustrate the sieve-stride optimization algorithm, Table 1 presents detailed pseudocode.

Table 1. Sieve Optimization Pseudocode

**Algorithm:** Sieve Step Optimization

**Input:** qs, gray_addi, corr
**Output:** Update sieve array and flags array in qs

```
1: Initialize sieve and flags to 0
2: Set end to sieve array end address
3: For i from qs->list[3] to qs->list[4]-1:
4:     Get prime, size, calculate offset co
5:     Update two root positions
6:     Set p0 = sieve + root_0, p1 = sieve + root_1
7:     For p0:
8:         bound = end - 4*prime
9:         While (p0 < bound) do
10:            Batch update p0[0] - p0[3*prime]
11:            p0 += 4*prime
12:        End While
13:        While (p0 < end) do
14:            Update *p0, p0 += prime
15:        End While
16:    Same operation for p1 (step7-step15)
17: End For
18: For i from qs->list[4] to qs->base.length-1:
19:    Get parameters; update root positions;
20:    Set p0, p1 (same as step3-6)
21:    For p0:
22:        bound = end - 4*prime
23:        While (p0 < bound) do
24:            Batch set flags, sieve values
25:            p0 += 4*prime
26:        End While
27:        While (p0 < end) do
28:            Set single flag; Update *p0
30:        End While
31:    Same operation for p1 (step21-30)
32: End For
```

On the other hand, the original program alternatively processed the two roots of the polynomial $f(x)$ $root_{1,2}$, which results in the problem of memory access skipping. This problem can be solved by processing $root_{1,2}$ separately. The optimization effect of the sieving method contributes the most the overall optimization. The optimization results are shown in Figure 2.

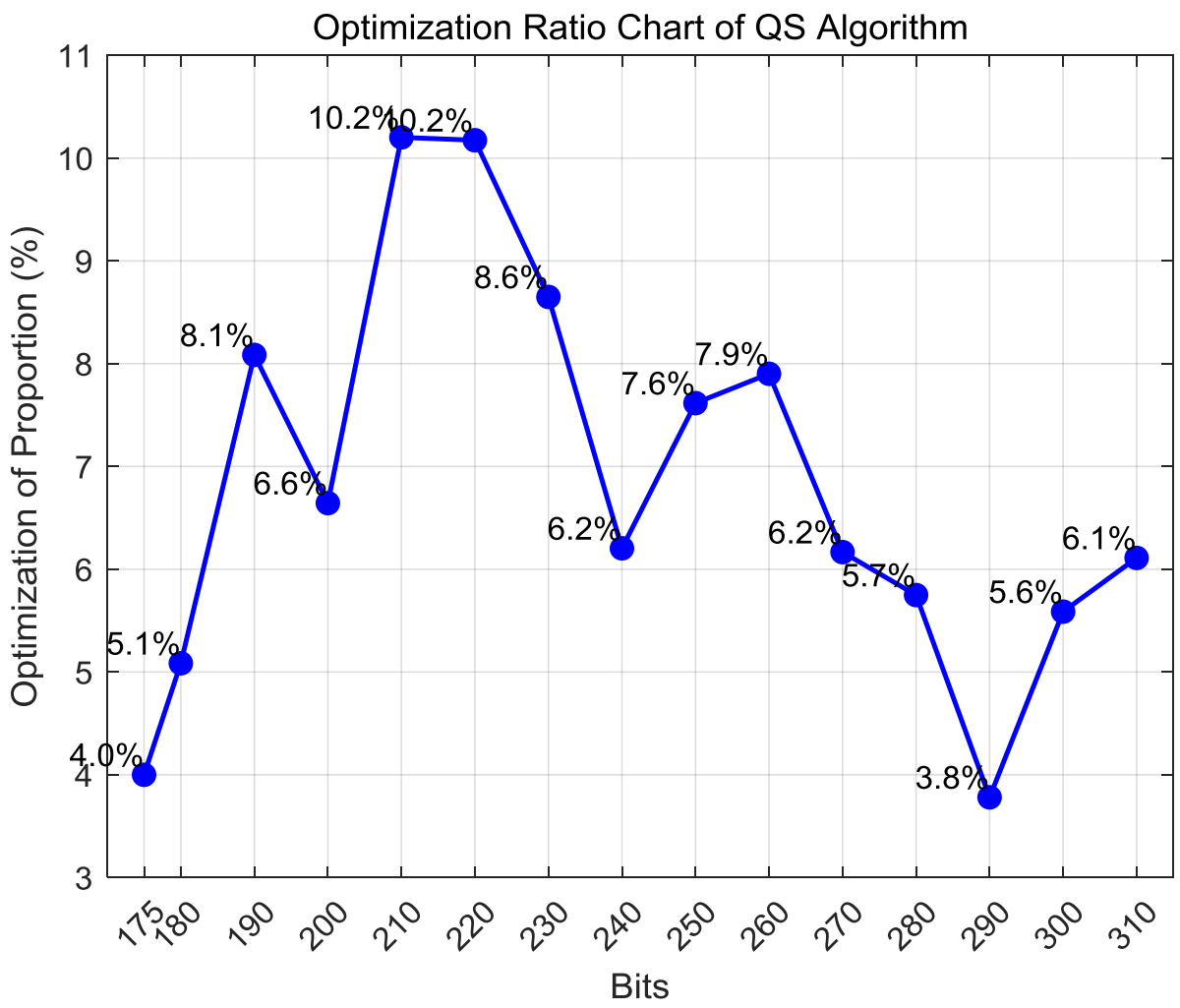

Fig. 2: Results of the optimization on the sieving method

Another optimization is implemented by parameter lookup tables, [14]. We replace the original traversal lookup with a binary search scheme. In order to determine the search range, we determine the range of the parameter table and find the last valid index by accumulating the flags. Applying binary search, our algorithm probes either the left or right half of the parameter table. This step precisely locates the target parameter interval, enabling linear interpolation within the identified range. Based on the above operations, we can construct a functional relationship between the bit-length of $kN$ and different parameters. Finally, we calibrate and align the output results.

After optimization, the time complexity of the loop can be reduced from $O(n)$ to $O(logn)$. At the same time, the boundary judgment of the parameter table is carried out to prevent overflow, which makes the project run more stably. This part has about 5% optimization effect of the overall improvement.

Table 2 presents the test results before and after optimization for large integer $N$ of different size.

Table 2. Decomposition time of $N$ with different digits

| **The number of digits of *N*(bit)** | **Original program running time (s)** | **Modify the sieve step size (s)** | **Modify parameter lookup and sieve step size (s)** |
|---|---|---|---|
| 70 | 0.03 | 0.03 | 0.05 |
| 80 | 0.05 | 0.04 | 0.02 |
| 90 | 0.03 | 0.03 | 0.04 |
| 100 | 0.03 | 0.03 | 0.02 |
| 110 | 0.03 | 0.02 | 0.03 |
| 120 | 0.02 | 0.02 | 0.06 |
| 130 | 0.07 | 0.06 | 0.09 |
| 140 | 0.09 | 0.09 | 0.11 |
| 150 | 0.10 | 0.11 | 0.40 |
| 160 | 0.48 | 0.47 | 0.75 |
| 170 | 0.73 | 0.70 | 1.05 |

Table 3. Decomposition time of N with different digits

| **The number of digits of N(bit)** | **Original code decomposition time (s)** | **Modify the sieve step size (s)** | **Modify parameter lookup and sieve step size (s)** |
|---|---|---|---|
| 150-170 | 8.27 | 8.01 | 7.38 |

Table 4. Decomposition time of N with different digits

| **The number of digits of N(bit)** | **Original code decomposition time (s)** | **Modify the sieve step size (s)** | **Modify parameter lookup and sieve step size (s)** |
|---|---|---|---|
| 175 | 1.00 | 0.96 | 0.95 |
| 180 | 1.18 | 1.12 | 1.10 |
| 190 | 2.35 | 2.16 | 1.93 |
| 200 | 5.72 | 5.34 | 4.55 |
| 210 | 8.92 | 8.01 | 7.72 |
| 220 | 25.85 | 23.22 | 22.78 |
| 230 | 41.39 | 37.81 | 35.38 |
| 240 | 76.09 | 71.37 | 65.00 |
| 250 | 153.36 | 141.68 | 144.49 |
| 260 | 271.47 | 250.02 | 251.24 |
| 270 | 453.42 | 425.46 | 407.05 |
| 280 | 1297.69 | 1223.09 | 1155.95 |
| 290 | 1931.72 | 1858.69 | 1841.29 |
| 300 | 12055.44 | 11381.88 | 11389.94 |
| 310 | 15404.66 | 14463.42 | 12475.35 |

For individual data, the optimization effect of data less than 170 bits is not obvious. So we take the total time consumed by factorization of different size $N$ from 150 bits to 170 bits. Table 3 presents the total factorization time consumed for 21 large integers ranging from 150 to 170 bits, decomposed under different optimization approaches. Then a significant efficiency improvement can be seen.

As a compiling detail, the $.exe$ files generated by executing " $\mathrm{gcc} - \mathrm{Wall} - \mathrm{pedantic} - \mathrm{O3} - \mathrm{std} = \mathrm{c99}$ C:\Users\hp\Desktop\C \main.c − o factor.exe " are faster than the $.exe$ files generated by "−O2" and "−O1". It's about 10%-15% faster, and all of the above data is used for .exe tests generated by "−O3".

For large numbers with more than 175 bits, a significant improvement in efficiency is observed; detailed data are presented in Table 4. But the efficiency improvement of different digits is different (5%-20% difference).

After modifying the parameter searching and sieving step size, the optimization ratio is shown in Figure 3.

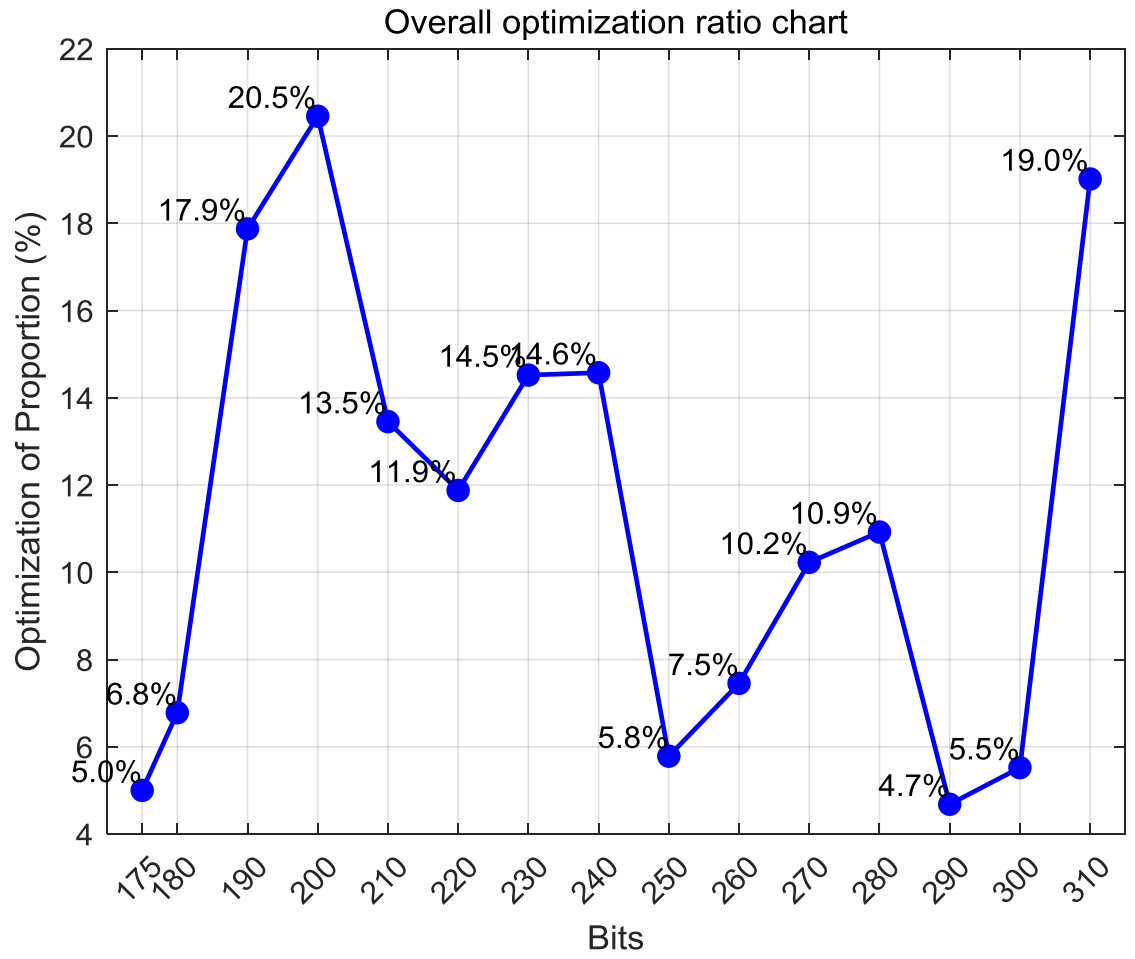

Fig. 3: The overall optimization effect of large integer factorization

## 5 Conclusion

Large integer factorization has important theoretical and practical significance in cryptography. It is closely related to the security of RSA algorithm. The method presented in this paper improves the efficiency of the existing quadratic sieve program by the optimization of sieving step size and parameter lookup tables. Future work will focus on optimizing algorithmic loop steps, potentially employing MPI-based multi-machine parallel computing to parallelize the most time-consuming loops. Additionally, the codebase can be ported to GPU architectures for enhanced performance.

**Contribution of Individual Authors to the Creation of a Scientific Article (Ghostwriting Policy)**

- Zihan Guan carried out the programming optimization and paper draft writing.
- Xiaodong Zhuang carried out the algorithm analysis.
- Nikos E. Mastorakis carried out the draft revising and proof reading.

**Sources of Funding for Research Presented in a Scientific Article or Scientific Article Itself**

No funding was received for conducting this study.

**Conflict of Interest**

The authors have no conflicts of interest to declare.

# APPENDIX

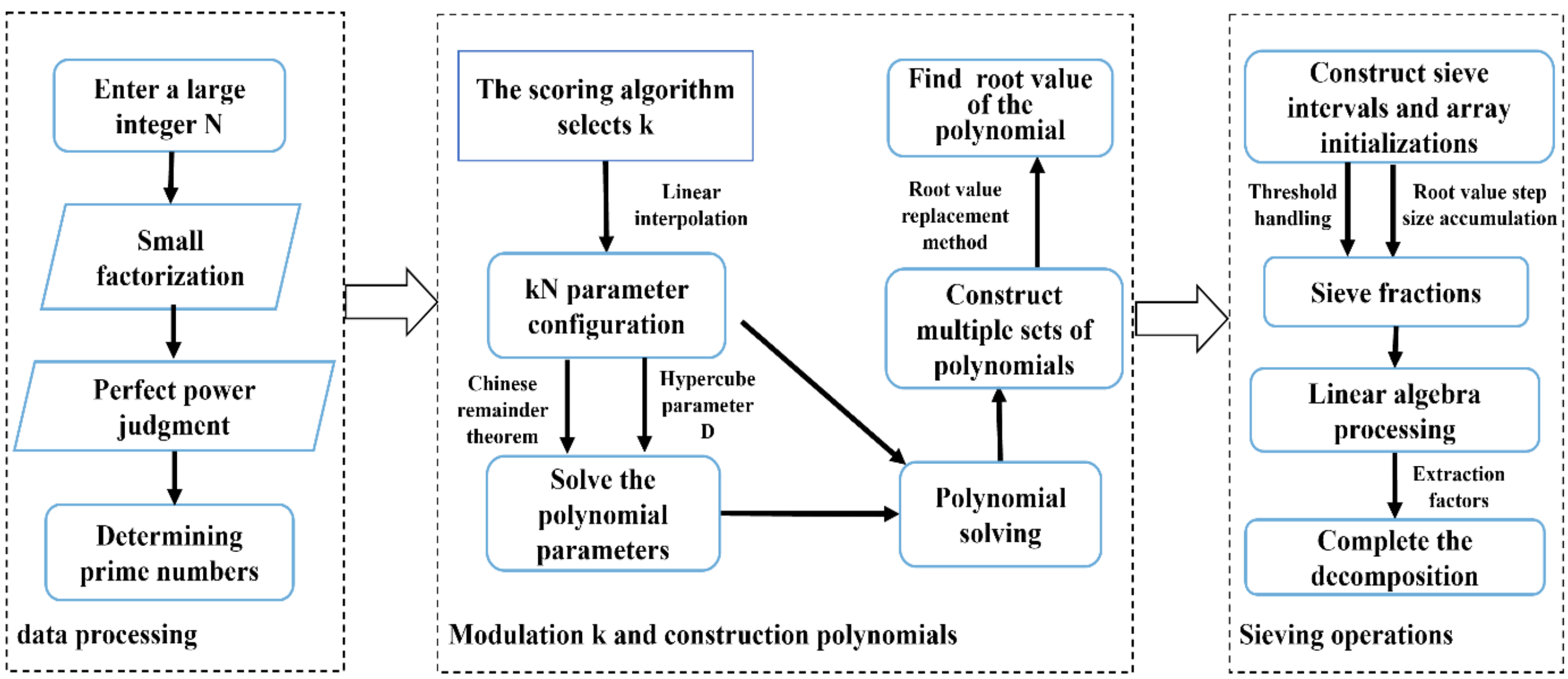

Fig. 1: The implementation process of the quadratic sieve algorithm